\documentclass[11pt]{amsart}
\input{diagrams}

\newtheorem{proposition}{Proposition}[section]
\newtheorem{lemma}[proposition]{Lemma}

\newtheorem{theorem}[proposition]{Theorem}

\theoremstyle{definition}

\theoremstyle{remark}

\newtheorem{remarks}[proposition]{Remarks}

\newcommand{\thlabel}[1]{\label{th:#1}}
\newcommand{\thref}[1]{Theorem~\ref{th:#1}}
\newcommand{\selabel}[1]{\label{se:#1}}
\newcommand{\seref}[1]{Section~\ref{se:#1}}
\newcommand{\lelabel}[1]{\label{le:#1}}
\newcommand{\leref}[1]{Lemma~\ref{le:#1}}
\newcommand{\prlabel}[1]{\label{pr:#1}}
\newcommand{\prref}[1]{Proposition~\ref{pr:#1}}

\newcommand{\relabel}[1]{\label{re:#1}}

\newcommand{\eqlabel}[1]{\label{eq:#1}}
\newcommand{\equref}[1]{(\ref{eq:#1})}

\newcommand{\Hom}{{\rm Hom}}

\newcommand{\End}{{\rm End}}

\newcommand{\im}{{\rm Im}\,}

\newcommand{\can}{{\rm can}}

\def\lan{\langle}
\def\ran{\rangle}
\def\ot{\otimes}

\newcommand{\Cc}{\mathcal{C}}
\newcommand{\Dd}{\mathcal{D}}

\newcommand{\Mm}{\mathcal{M}}

\def\text#1{{\rm {\rm #1}}}

\def\ol{\overline}
\def\ul{\underline}

\begin{document}
\title[Galois theory for weak Hopf algebras]{Galois theory for weak Hopf algebras}
\author{S. Caenepeel}
\address{Faculty of Engineering Sciences,
Vrije Universiteit Brussel, VUB, B-1050 Brussels, Belgium}
\email{scaenepe@vub.ac.be}
\urladdr{http://homepages.vub.ac.be/\~{}scaenepe/}
\author{E. De Groot}
\address{Faculty of Engineering Sciences,
Vrije Universiteit Brussel, VUB, B-1050 Brussels, Belgium}
\email{edegroot@vub.ac.be}
\urladdr{http://homepages.vub.ac.be/\~{}edegroot/}
\begin{abstract}
We develop Hopf-Galois theory for weak Hopf algebras, and recover
analogs of classical results for Hopf algebras. Our methods are based on
the recently introduced Galois theory for corings. We focus on the situatation
where the weak Hopf algebra is a groupoid algebra or its dual. We obtain
weak versions of strongly graded rings, and classical Galois extensions.
\end{abstract}

\maketitle

\section*{Introduction}
Weak bialgebras and Hopf algebras are generalizations of ordinary bialgebras
and Hopf algebras in the following sense: the defining axioms are the same,
but the multiplicativity of the counit and comultiplicativity of the unit
are replaced by weaker axioms. Perhaps the easiest example of a weak Hopf algebra
is a groupoid algebra; other examples are face algebras \cite{Hayashi},
quantum groupoids \cite{Ocneanu} and generalized Kac algebras \cite{Yamanouchi}.
A purely algebraic approach can be found in \cite{BohmNSI} and \cite{BohmS}.\\
The aim of this note is to develop Galois theory for weak Hopf algebras.
A possible strategy could be to try to adapt the methods from classical
Hopf-Galois theory to the weak situation. This would be unnecessarily complicated,
since a more powerful tool became available recently.\\
Corings and comodules over corings were introduced by Sweedler in \cite{S3}.
The notion was revived recently by Brzezi\'nski \cite {B2}, who saw the
importance of a remark made by Takeuchi, that Hopf modules and many of their
generalizations can be viewed as examples of comodules over corings. Many
applications came out of this idea, an overview is given in \cite{BW}.\\
One of the beautiful applications is a reformulation of descent and Galois
theory using the language of corings. This was initiated in \cite{B2},
and continued by the first author \cite{Caenepeel03} and Wisbauer \cite{W3}.
An observation in \cite{B2} is that weak Hopf modules can be viewed as
comodules over a coring, and this implies that the general theory of
Galois corings can be applied to weak Hopf algebras. This is what we will
do in \seref{2}.\\ 
In \seref{3}, we look at the special case where the 
weak Hopf algebra is a groupoid algebra. This leads us to the notions of
groupoid graded algebras and modules, and the following generalization of a result
of Ulbrich: a groupoid graded algebra is Galois if and only if it is strongly
graded.\\
In \seref{4}, we look at the dual situation, where the weak Hopf algebra is
the dual of a finite groupoid algebra. Then we have to deal with finite groupoid
actions, as we have to deal with finite group actions in classical Galois theory.
The results in this Section are similar to the ones obtained in \cite{CDG2},
where partial Galois theory (see \cite{DFP}) is studied using the language of corings.
Both are generalizations of classical Galois theory. In partial Galois theory, we
look at group actions on algebra $A$
that are not everywhere defined, and in weak Galois theory,
we look at groupoid actions. In both cases, we have a coring that is a direct factor
of $|G|$ copies of $A$, and the left dual of the coring is a Frobenius extension of $A$.

\section{Preliminary results}\selabel{1}
\subsection{Galois corings}\selabel{1.1}
Let $A$ be a ring. An $A$-coring $\Cc$ is a coalgebra in the category
${}_A\Mm_A$ of $A$-bimodules. Thus an $A$-coring
is a triple $\Cc=(\Cc,\Delta_\Cc,\varepsilon_\Cc)$, where
 $\Cc$ is an $A$-bimodule, and
$\Delta_\Cc:\ \Cc\to \Cc\ot_A\Cc$  and
 $\varepsilon_\Cc:\ \Cc\to A$ are $A$-bimodule maps
such that
\begin{equation}\eqlabel{1.1.1}
(\Delta_\Cc\ot_A \Cc)\circ \Delta_\Cc=
(\Cc\ot_A\Delta_\Cc )\circ \Delta_\Cc,
\end{equation}
and
\begin{equation}\eqlabel{1.1.2}
({\Cc}\ot_A\varepsilon_{\Cc})\circ \Delta_\Cc=
(\varepsilon_{\Cc}\ot_A{\Cc})\circ \Delta_\Cc=\Cc.
\end{equation}
We use the Sweedler-Heyneman notation for the comultiplication:
$$\Delta_\Cc(c)=c_{(1)}\ot_A c_{(2)}.$$
A right $\Cc$-comodule $M=(M,\rho )$ consists of a right $A$-module
$M$ together with a right $A$-linear map $\rho :\ M\to M\ot_A\Cc$
such that:
\begin{equation}\eqlabel{1.1.3}
(\rho \ot_A \Cc)\circ \rho =(M\ot_A\Delta_\Cc)\circ \rho, 
\end{equation}
and
\begin{equation}\eqlabel{1.1.4}
(M\ot_A\varepsilon_\Cc)\circ \rho =M.
\end{equation}
We then say that $\Cc$ coacts from the right on $M$, and we denote
$$\rho (m)=m_{[0]}\ot_A m_{[1]}.$$
A right $A$-linear map
$f:\ M\to N$ between two right $\Cc$-comodules $M$ and $N$ is called
right $\Cc$-colinear if $\rho (f(m))=f(m_{[0]})\ot m_{[1]}$,
for all $m\in M$. The category of right $\Cc$-comodules and $\Cc$-colinear
maps is denoted by $\Mm^\Cc$.\\
$x\in \Cc$ is called grouplike if $\Delta_\Cc(x)=x\ot x$ and
$\varepsilon_\Cc(x)=1$. Grouplike elements of $\Cc$ correspond bijectively
to right $\Cc$-coactions on $A$: if $A$ is grouplike, then we have the following
right $\Cc$-coaction $\rho$ on $A$: $\rho(a)=xa$.\\
Let $(\Cc,x)$ be a coring with a fixed grouplike element. For $M\in \Mm^{\Cc}$, we call
$$M^{{\rm co}\Cc}=\{m\in M~|~\rho (m)=m\ot_A x\}$$
the submodule of coinvariants of $M$. Observe that
$$A^{{\rm co}\Cc}=\{b\in A~|~bx=xb\}$$
is a subring of $A$. Let $i:\ B\to A$ be a ring morphism.
$i$ factorizes through $A^{{\rm co}\Cc}$ if and only if
$$x\in G(\Cc)^B=\{x\in G(\Cc)~|~xb=bx,~{\rm for~all~}b\in B\}.$$ 
We then have
a pair of adjoint functors
$(F,G)$, respectively between the categories $\Mm_B$ and
$\Mm^\Cc$ and the categories ${}_B\Mm$ and
${}^\Cc\Mm$. For 
$N\in \Mm_B$ and $M\in \Mm^\Cc$,
$$F(N)=N\ot_B A~~{\rm and}~~G(M)=M^{{\rm co}\Cc}.$$
The unit and counit of the adjunction are
$$\nu_N:\ N\to (N\ot_BA)^{{\rm co}\Cc},~~\nu_N(n)=n\ot_B1;$$
$$\zeta_M:\ M^{{\rm co}\Cc}\ot_B A\to M,~~
\zeta_M(m\ot_B a)=ma.$$
We have a similar pair of adjoint functors $(F',G')$ between the categories of left 
$B$-modules and left $\Cc$-comodules.\\
Let $i:\ B\to A$ be a morphism of rings. The associated canonical coring is
$\Dd=A\ot_BA$, with comultiplication and counit given by the formulas
$$\Delta_\Dd:\ \Dd\to \Dd\ot_A\Dd\cong A\ot_B A\ot_B A,~~
\Delta_\Dd(a\ot_B a')=a\ot_B1\ot_Ba'$$
and
$$\varepsilon_\Dd:\ \Dd=A\ot_B A\to A,~~\varepsilon_\Dd(a\ot_Ba')=aa'.$$
If $i:\ B\to A$ is pure as a morphism of left and right $B$-modules, then
the categories $\Mm_B$ and $\Mm^\Dd$ are equivalent.\\
Let $(\Cc,x)$ be a coring with a fixed grouplike element, and $i:\ B\to A^{{\rm co}\Cc}$
a ring morphism. We then have a morphism of corings
$${\rm can}:\ \Dd=A\ot_BA\to \Cc,~~{\rm can}(a\ot_B a')=axa'.$$
If $F$ is fully faithful, then $B\cong A^{{\rm co}\Cc}$; if $G$ is fully
faithful, then ${\rm can}$ is  an isomorphism.
$(\Cc,x)$ is called a Galois coring if  
${\rm can}:\ A\ot_{A^{{\rm co}\Cc}}A\to \Cc$ is
bijective. From \cite{Caenepeel03}, we recall the following results.

\begin{theorem}\thlabel{1.2}
Let $(\Cc,x)$ be an $A$-coring with fixed grouplike element, and
$T=A^{{\rm co}\Cc}$. Then the following statements are equivalent.
\begin{enumerate}
\item $(\Cc,x)$ is Galois and $A$ is faithfully flat as a left $T$-module;
\item $(F,G)$ is an equivalence
and $A$ is flat as a left $T$-module.
\end{enumerate}
\end{theorem}

Let $(\Cc,x)$ be a coring with a fixed grouplike element, and take
$T=A^{{\rm co}\Cc}$. Then ${}^*\Cc={}_A\Hom(\Cc,A)$ is a ring, with
multiplication given by
\begin{equation}\eqlabel{1.2.1a}
(f\# g)(c)=g(c_{(1)}f(c_{(2)})).
\end{equation}
We have a morphism of rings $j:\ A\to {}^*\Cc$, given by
$$j(a)(c)=\varepsilon_\Cc(c)a.$$
This makes ${}^*\Cc$ into an $A$-bimodule,
via the formula
$$(afb)(c)=f(ca)b.$$
Consider the left dual of the canonical map:
$${}^*{\rm can}:\ {}^*\Cc\to {}^*\Dd\cong {}_T\End(A)^{\rm op},~~
{}^*{\rm can}(f)(a)=f(xa).$$
We then have the following result.

\begin{proposition}\prlabel{1.3}
If $(\Cc,x)$ is Galois, then ${}^*{\rm can}$ is an isomorphism.
The converse property holds if $\Cc$ and $A$ are finitely generated
projective, respectively as a left $A$-module, and a left $T$-module.
\end{proposition}

Let $Q=\{q\in {}^*\Cc~|~c_{(1)}q(c_{(2)})=q(c)x,~{\rm for~all~}c\in \Cc\}$.
A straightforward computation shows that $Q$ is a $({}^*\Cc,T)$-bimodule.
Also $A$ is a left $(T,{}^*\Cc)$-bimodule; the right ${}^*\Cc$-action is
induced by the right $\Cc$-coaction: $a\cdot f=f(xa)$. Now
consider the maps
\begin{eqnarray}
\tau:&& A\ot_{{}^*\Cc} Q\to T,~~\tau(a\ot_{{}^*\Cc} q)=q(xa);\eqlabel{1.3.1}\\
\mu:&& Q\ot_T A\to {}^*\Cc,~~\mu(q\ot_T a)=q\#i(a).\eqlabel{1.3.2}
\end{eqnarray}
With this notation, we have the following property (see \cite{CVW}).

\begin{proposition}\prlabel{1.4}
$(T,{}^*\Cc,A,Q,\tau,\mu)$ is a Morita context. The map $\tau$ is surjective
if and only if there exists $q\in Q$ such that $q(x)=1$.
\end{proposition}

We also have (see \cite{Caenepeel03}):

\begin{theorem}\thlabel{1.5}
Let $(\Cc,x)$ be a coring with fixed grouplike element, and assume
that $\Cc$ is a left $A$-progenerator. We take a subring
$B$ of $T=A^{{\rm co}\Cc}$, and consider the map
$${\rm can}:\ \Dd=A\ot_{B}A\to \Cc,~~{\rm can}(a\ot_{T}a')=axa'$$
Then the following statements are equivalent:
\begin{enumerate}
\item \begin{itemize}
\item ${\rm can}$ is an isomorphism;
\item $A$ is faithfully flat as a left $B$-module.
\end{itemize}
\item \begin{itemize}
\item ${}^*{\rm can}$ is an isomorphism;
\item $A$ is a left $B$-progenerator.
\end{itemize}
\item \begin{itemize}
\item $B=T$;
\item the Morita context $(B,{}^*\Cc,A,Q,\tau,\mu)$ is strict.
\end{itemize}
\item \begin{itemize}
\item $B=T$;
\item $(F,G)$ is an equivalence of categories.
\end{itemize}
\end{enumerate}
\end{theorem}

\subsection{Weak Hopf algebras}\selabel{1.2}
Let $k$ be a commutative ring. Recall that a weak $k$-bialgebra is a $k$-module
with a $k$-algebra structure $(\mu,\eta)$ and a $k$-coalgebra structure
$(\Delta,\varepsilon)$ such that
$$\Delta(hk)=\Delta(h)\Delta(k),$$
for all $h,k\in H$, and
\begin{eqnarray}
\Delta^2(1)&=& 1_{(1)}\ot 1_{(2)}1_{(1')}\ot 1_{(2')}=
1_{(1)}\ot 1_{(1')}1_{(2)}\ot 1_{(2')},\eqlabel{1.6.1}\\
\varepsilon(hkl)&=& \varepsilon(hk_{(1)})\varepsilon(k_{(2)}l)=
\varepsilon(hk_{(2)})\varepsilon(k_{(1)}l),\eqlabel{1.6.2}
\end{eqnarray}
for all $h,k,l\in H$. We use the Sweedler-Heyneman notation for the
comultiplication, namely
$$\Delta(h)=h_{(1)}\ot h_{(2)}=h_{(1')}\ot h_{(2')}.$$
If $H$ is a weak bialgebra, then we have idempotent maps
$\Pi^L,\Pi^R,\ol{\Pi}^L,\ol{\Pi}^R:\ H\to H$ given by
\begin{eqnarray*}
\Pi^L(h)&=& \varepsilon(1_{(1)}h)1_{(2)}\\
\Pi^R(h)&=& 1_{(1)}\varepsilon(h1_{(2)})\\
\ol{\Pi}^L(h)&=&1_{(1)}\varepsilon(1_{(2)}h)\\
\ol{\Pi}^R(h)&=&\varepsilon(h1_{(1)})1_{(2)}
\end{eqnarray*}
We have that $H^L=\im (\Pi^L)=\im (\ol{\Pi}^R)$ and
$H^R=\im (\Pi^R)=\im (\ol{\Pi}^L)$.\\
Let $H$ be a weak bialgebra. Recall from \cite{Bohm} and \cite{CDG}
that a right $H$-comodule algebra $A$ is a $k$-algebra with a
right $H$-comodule structure $\rho$ such that $\rho(a)\rho(b)=\rho(ab)$,
for all $a,b\in A$, and such that
the following equivalent statements hold (see \cite[Prop. 4.10]{CDG}):
\begin{eqnarray}
&&\rho^2(1)= \sum 1_{[0]}\ot 1_{[1]}1_{(1)}\ot 1_{(2)}\eqlabel{1.6.3}\\
&&\rho^2(1)= \sum 1_{[0]}\ot 1_{(1)}1_{[1]}\ot 1_{(2)}\eqlabel{1.6.4}\\
&&\sum a_{[0]}\ot\ol{\Pi}^R(a_{[1]})=\sum a1_{[0]}\ot 1_{[1]}\eqlabel{1.6.5}\\
&&\sum a_{[0]}\ot{\Pi}^L(a_{[1]})=\sum 1_{[0]}a\ot 1_{[1]}\eqlabel{1.6.6}\\
&&\sum 1_{[0]}\ot\ol{\Pi}^R(1_{[1]})=\rho^r(1)\eqlabel{1.6.7}\\
&&\sum 1_{[0]}\ot{\Pi}^L(1_{[1]})=\rho^r(1)\eqlabel{1.6.8}\\
&&\rho^r(1)\in A\ot H^L\eqlabel{1.6.9}
\end{eqnarray}

\begin{lemma}\lelabel{1.7}
Let $H$ be a weak bialgebra, and $A$ a right $H$-comodule algebra.
Then
\begin{eqnarray}
\varepsilon(h_{(1)}1_{[1]})1_{[0]}\ot h_{(2)}&=&
1_{[0]}\ot h 1_{[1]}\eqlabel{1.7.1}\\
\varepsilon(h1_{[1]})1_{[0]}a&=& \varepsilon(ha_{[1]})a_{[0]}\eqlabel{1.7.2}
\end{eqnarray}
for all $h\in H$ and $a\in A$.
\end{lemma}

\begin{proof}
This is a special case of \cite[Theorem 4.14]{CDG}. Details are as follows.
\begin{eqnarray*}
&&\hspace*{-2cm} \varepsilon(h_{(1)}1_{[1]})1_{[0]}\ot h_{(2)}=
\varepsilon(h_{(1)}1_{(1)}1_{[1]})1_{[0]}\ot h_{(2)}1_{(2)}\\
&=& \varepsilon(h_{(1)}1_{[1]})1_{[0]}\ot h_{(2)}1_{[2]}=1_{[0]}\ot h 1_{[1]};\\
&&\hspace*{-2cm}\varepsilon(h1_{[1]})1_{[0]}a=
\varepsilon(h1_{[2]})\varepsilon(1_{[1]}a_{[1]})1_{[0]}a_{[0]}\\
&=& \varepsilon(h1_{[1]}a_{[1]})1_{[0]}a_{[0]}=\varepsilon(ha_{[1]})a_{[0]}.
\end{eqnarray*}
\end{proof}

A weak Hopf algebra is a weak bialgebra together with a map $S:\ H\to H$,
called the antipode, satisfying the following conditions, for all $h\in H$:
\begin{equation}\eqlabel{1.6.10}
h_{(1)}S(h_{(2)})=\Pi^L(h);~~S(h_{(1)})h_{(2)}=\Pi^R(h);
\end{equation}
\begin{equation}\eqlabel{1.6.11}
S(h_{(1)})h_{(2)}S(h_{(3)})=S(h).
\end{equation}

\section{Weak Hopf-Galois extensions}\selabel{2}
Let $H$ be a weak bialgebra, and $A$ a right $H$-comodule algebra.
We then have a projection
$$g:\ A\ot H\to A\ot H,~~g(a\ot h)=a 1_{[0]}\ot h 1_{[1]}.$$
Observe that, for any $a\in A$,
$$\rho(a)=a_{[0]}\ot a_{[1]}=a_{[0]}1_{[0]}\ot a_{[1]}1_{[1]}=g(\rho(a))
\in \im(g).$$

\begin{lemma}\lelabel{2.1}
Let $H$ be a weak bialgebra. Then $\Cc=\im(g)$ is an $A$-coring, with structure
maps
\begin{eqnarray*}
b'(a 1_{[0]}\ot h 1_{[1]})b&=& b'ab_{[0]}\ot h b_{[1]}\eqlabel{2.1.1};\\
\Delta_\Cc(1_{[0]}\ot h 1_{[1]})&=&
(1_{[0]}\ot h_{(1)} 1_{[1]})\ot_A (1\ot h_{(2)} )\eqlabel{2.1.2};\\
\varepsilon_\Cc(1_{[0]}\ot h 1_{[1]})&=&1_{[0]}\varepsilon(h1_{[1]})\eqlabel{2.1.3}.
\end{eqnarray*}
$\rho(1)=1_{[0]}\ot 1_{[1]}$ is a grouplike element in $\Cc$.
\end{lemma}

\begin{proof}
This is a special case of a result in \cite[Prop. 2.3]{B1}. The fact that
$\Delta_\Cc(1_{[0]}\ot h 1_{[1]})\in \Cc\ot_A\Cc$ follows from
\begin{equation}\eqlabel{2.1.4}
\Delta_\Cc(1_{[0]}\ot h 1_{[1]})=
(1_{[0]}\ot h_{(1)}1_{[1]})\ot_A ((1_{[0']}\ot h_{(2)}1_{[1']}).
\end{equation}
Indeed,
\begin{eqnarray*}
&&\hspace*{-2cm}
(1_{[0]}\ot h_{(1)}1_{[1]})\ot_A ((1_{[0']}\ot h_{(2)}1_{[1']})\\
&=&(1_{[0]}1_{[0']}\ot h_{(1)}1_{[1]}1_{[1']})\ot_A (1\ot h_{(2)}1_{[2']})\\
&=& (1_{[0]}1_{[0']}\ot h_{(1)}1_{[1]}1_{[1']}1_{(1)})\ot_A (1\ot h_{(2)}1_{(2)})\\
&=&(1_{[0]}\ot h_{(1)}1_{[1]}1_{(1)})\ot_A (1\ot h_{(2)}1_{(2)})\\
&=& (1_{[0]}\ot h_{(1)}1_{(1)}1_{[1]})\ot_A (1\ot h_{(2)}1_{(2)})=\rho(a).
\end{eqnarray*}
Taking $h=1$ in \equref{2.1.3} and \equref{2.1.4}, we see that $\varepsilon_\Cc(\rho(1))=1$ and $\Delta_\Cc(\rho(1))=\rho(1)\ot_A \rho(1)$.
\end{proof}

We call $(M,\rho)$ a relative right $(A,H)$-Hopf module if $M$ is a right $A$-module,
$(M,\rho)$ is a right $H$-comodule, and
\begin{equation}\eqlabel{2.2.1}
\rho(ma)=m_{[0]}a_{[0]}\ot m_{[1]}a_{[1]},
\end{equation}
for all $m\in M$ and $a\in A$. $\Mm_A^H$ is the category of relative Hopf modules
and right $A$-linear $H$-colinear maps.

The following result is a special case of  \cite[Prop. 2.3]{B1}. Since it is
essential in what follows, we give a sketch of proof.

\begin{proposition}\prlabel{2.2}
The category $\Mm_A^H$ of relative Hopf modules is isomorphic to the
category $\Mm^\Cc$ of right $\Cc$-comodules over the coring $\Cc=\im g$.
\end{proposition}

\begin{proof}
Take $(M,\tilde{\rho})\in \Mm^\Cc$. Let $\rho$ be the composition
$$\rho:\ M\rTo^{\tilde{\rho}} M\ot_A\Cc\rTo^{} M\ot_A(A\ot H)
\rTo^{\cong} M\ot H.$$
If
$$\tilde{\rho}(m)=\sum_i m_i \ot_A a_i 1_{[0]}\ot h_i1_{[1]},$$
then
$$\rho(m)=m_{[0]}\ot m_{[1]}= \sum_i m_i a_i 1_{[0]}\ot h_i1_{[1]}.$$
Then $\tilde{\rho}$ is determined by $\rho$, since
\begin{eqnarray}
\tilde{\rho}(m)&=&\sum_i m_i \ot_A a_i 1_{[0]}\ot h_i1_{[1]}\nonumber\\
&=& \sum_i m_i \ot_A a_i 1_{[0]}1_{[0']}\ot h_i1_{[1]}1_{[0']}\nonumber\\
&=& \sum_i m_i  a_i 1_{[0]}\ot_A 1_{[0']}\ot h_i1_{[1]}1_{[0']}\nonumber\\
&=& m_{[0]}\ot_A 1_{[0]}\ot m_{[1]}1_{[1]}.\eqlabel{2.2.2}
\end{eqnarray}
From the fact that $\tilde{\rho}$ is right $A$-linear, it follows that
\begin{eqnarray*}
&&\hspace*{-2cm}
\tilde{\rho}(ma)=\tilde{\rho}(m)a=
m_{[0]}\ot_A (1_{[0]}\ot m_{[1]}1_{[1]})a\\
&=& m_{[0]}\ot_A a_{[0]}\ot m_{[1]}a_{[1]},
\end{eqnarray*}
so
$$\rho(ma)=m_{[0]} a_{[0]}\ot m_{[1]}a_{[1]}.$$
From the fact that $(M\ot_A\varepsilon_\Cc)(\tilde{\rho}(m))=m$ and
\equref{2.2.2}, it follows that
$$m=m_{[0]} 1_{[0]}\ot \varepsilon(m_{[1]}1_{[1]})
=m_{[0]} \ot \varepsilon(m_{[1]}).$$
For all $m\in M$, we have the following equality in $M\ot_A(A\ot H)\ot_A(A\ot H)$;
\begin{equation}\eqlabel{2.2.3}
(M\ot_A\Delta_{\Cc})(\tilde{\rho}(m))=(\tilde{\rho}\ot \Cc)(\tilde{\rho}(m))
\end{equation}
The left hand side of \equref{2.2.3} is
$$m_{[0][0]}\ot_A(1_{[0']}\ot m_{[0][1]}1_{[1']})\ot_A
(1_{[0]}\ot m_{[1]}1_{[1]}).$$
The image in $M\ot_A(A\ot H\ot H)$ is
\begin{eqnarray*}
&&\hspace*{-2cm}
m_{[0][0]}\ot_A(1_{[0']}1_{[0]}\ot m_{[0][1]}1_{[1']}1_{[1]}1_{(1)}\ot 
m_{[1]}1_{(2)})\\
&=& m_{[0][0]}\ot_A(1_{[0]}\ot m_{[0][1]}1_{[1]}1_{(1)}\ot 
m_{[1]}1_{(2)}),
\end{eqnarray*}
and in $M\ot H\ot H$:
\begin{eqnarray*}
&&\hspace*{-2cm}
m_{[0][0]}1_{[0]}\ot m_{[0][1]}1_{[1]}1_{(1)}\ot m_{[1]}1_{(2)}\\
&=& m_{[0][0]}1_{[0]}\ot m_{[0][1]}1_{[1]}\ot m_{[1]}1_{[2]}\\
&=& m_{[0][0]}1_{[0][0]}\ot m_{[0][1]}1_{[0][1]}\ot m_{[1]}1_{[1]}\\
&=& m_{[0][0]}\ot m_{[0][1]}\ot m_{[1]}.
\end{eqnarray*}
The right hand side of \equref{2.2.3} is
$$m_{[0]}\ot_A (1_{[0]}\ot m_{[1](1)}1_{[1]})\ot_A (1_{[0']}\ot m_{[1](2)}1_{[1']}).$$
The image in $M\ot_A(A\ot H\ot H)$ is
\begin{eqnarray*}
&&\hspace*{-2cm}
m_{[0]}\ot_A (1_{[0]}1_{[0']}\ot m_{[1](1)}1_{[1]}1_{[1']}1_{(1)}\ot m_{[1](2)}1_{(2)})\\
&=& m_{[0]}\ot_A (1_{[0]}\ot m_{[1](1)}1_{[1]}1_{(1)}\ot m_{[1](2)}1_{(2)}),
\end{eqnarray*}
and in $M\ot H\ot H$:
\begin{eqnarray*}
&&\hspace*{-2cm}
m_{[0]}1_{[0]}\ot m_{[1](1)}1_{[1]}1_{(1)}\ot m_{[1](2)}1_{(2)}\\
&=&m_{[0]}1_{[0]}\ot m_{[1](1)}1_{[1]}\ot m_{[1](2)}1_{[2]}\\
&=&m_{[0]}1_{[0]}\ot m_{[1](1)}1_{[1](1)}\ot m_{[1](2)}1_{[1](2)}\\
&=& m_{[0]}\ot m_{[1](1)}\ot m_{[1](2)}.
\end{eqnarray*}
It follows that
$$m_{[0][0]}\ot m_{[0][1]}\ot m_{[1]}=m_{[0]}\ot m_{[1](1)}\ot m_{[1](2)},$$
and $\rho$ is coassociative. Conversely, given a relative Hopf module $(M,\rho)$,
we define $\tilde{\rho}:\ M\to M\ot_A\Cc$ using \equref{2.2.2}. Straightforward
computations show that $(M,\tilde{\rho})\in\Mm^\Cc$.
\end{proof}

Let $(M,\rho)$ be a relative Hopf module, and $(M,\tilde{\rho})$ the
corresponding $\Cc$-comodule. Then $m\in M^{{\rm co}\Cc}$ if and only if
$$\tilde{\rho}(m)=m\ot_A1_{[0]}\ot 1_{[1]}$$
if and only if
$$\rho(m)=m1_{[0]}\ot 1_{[1]}.$$
We conclude that $M^{{\rm co}\Cc}=M^{{\rm co}H}$, which is by definition
$$M^{{\rm co}H}=\{m\in M~|~\rho(m)=m1_{[0]}\ot 1_{[1]}\}.$$
The grouplike element $\rho(1)$ induces a right $\Cc$-coaction $\tilde{\rho}$
on $A$:
$$\tilde{\rho}(a)=1\ot_A (1_{[0]}\ot 1_{[1]})a= 1\ot_A (a_{[0]}\ot a_{[1]}).$$
The corresponding $H$-coaction on $A$ is the original $\rho$:
$$\rho(a)=a_{[0]}\ot a_{[1]}.$$
Observe that
$$T=A^{{\rm co}H}=\{a\in A~|~a1_{[0]}\ot 1_{[1]}\}.$$
Let $i:\ B\to T$ be a ring morphism. We have seen in \seref{1.1} that we have
a pair of adjoint functors $(F,G)$:
$$F:\ \Mm_B\to \Mm^H_A,~~F(N)=N\ot_B A;$$
$$G:\ \Mm^H_A\to \Mm_B,~~G(N)=N^{{\rm co}H}.$$
$F(N)=N\ot_B A$ is a relative Hopf module via
$$\rho(n\ot_B a)=n\ot_B a_{[0]}\ot a_{[1]}.$$
The canonical map $\can:\ A\ot_B A\to \Cc$ takes the form
\begin{equation}\eqlabel{2.3.1}
\can(a\ot_B b)=a(1_{[0]}\ot 1_{[1]})b=ab_{[0]}\ot b_{[1]}.
\end{equation}

We call $A$ a weak $H$-Galois extension of $A^{{\rm co}H}$ if
$(\Cc,\rho(1))$ is a Galois coring, that is, $\can:\
A\ot_{A^{{\rm co}H}} A\to \Cc$ is an isomorphism. From \thref{1.2},
we immediately have the following result.

\begin{proposition}\prlabel{2.3a}
Let $H$ be a weak bialgebra, and $A$ a right $H$-comodule algebra. 
Then the following assertions are equivalent
\begin{enumerate}
\item $A$ is $A$ a weak $H$-Galois extension of $T$ and faithfully flat
as a left $B$-module;
\item $(F,G)$ is an equivalence and $A$ is flat as a left $T$-module.
\end{enumerate}
\end{proposition}

Our next aim is to compute ${}^*\Cc$.

\begin{proposition}\prlabel{2.3}
Let $H$ be a weak bialgebra, $A$ a right $H$-comodule algebra, and
$\Cc=\im(g)$ the $A$-coring that we introduced above. Then
\begin{eqnarray*}
&&\hspace*{-2cm}
{}^*\Cc={}_A\Hom(\Cc,A)\cong \ul{\Hom}(H,A)\\
\{f\in \Hom(H,A)~|~f(h)=:1_{[0]}f(h1_{[1]}),~{\rm for~all~}h\in H\},
\end{eqnarray*}
with multiplication rule
\begin{equation}\eqlabel{2.3.2}
(f\# g)(h)=f(h_{(2)})_{[0]}g(h_{(1)}f(h_{(2)})_{[1]}).
\end{equation}
\end{proposition}

\begin{proof}
We define $\alpha:\ {}^*\Cc\to \ul{\Hom}(H,A)$ by
$$\alpha(\varphi)(h)=\varphi(1_{[0]}\ot h1_{[1]}).$$
$\alpha(\varphi)=f\in \ul{\Hom}(H,A)$ since
\begin{eqnarray*}
0&=& \alpha(\varphi)((1\ot h)(\rho(1)-\rho(1)^2))\\
&=&\alpha(\varphi)(1_{[0]}\ot h1_{[1]})
-\alpha(\varphi)(1_{[0]}1_{[0']}\ot h1_{[1]}1_{[1']})\\
&=& f(h)-1_{[0]}f(h1_{[1]}).
\end{eqnarray*}
Now we define $\beta:\ \ul{\Hom}(H,A)\to {}^*\Cc$ by
$$\beta(f)(a1_{[0]}\ot h1_{[1]})=af(h).$$
We have to show that $\beta$ is well-defined. Assume that
$\sum_i a_i1_{[0]}\ot h_i1_{[1]}= 0$
in $A\ot H$. Then
$$\sum_i a_i\ot h_i=
\sum_i a_i\ot h_i-\sum_i a_i1_{[0]}\ot h_i1_{[1]},$$
hence
$$\sum_i a_if(h_i)=\sum_i a_if(h_i)- \sum_i a_i1_{[0]}f( h_i1_{[1]})=0.$$
A straightforward verification shows that $\alpha$ is the inverse of $\beta$.
Using $\alpha$ and $\beta$, we can transport the multiplication on ${}^*\Cc$
to $\ul{\Hom}(H,A)$. This gives
\begin{eqnarray*}
&&\hspace*{-2cm}
(f\# g)(h)=(\varphi\#\psi)(1_{[0]}\ot h1_{[1]})\\
&=& \psi((1_{[0]}\ot h_{(1)}1_{[1]})\varphi(1_{[0']}\ot h_{(2)}1_{[1']})\\
&=& \psi((1_{[0]}\ot h_{(1)}1_{[1]})f(h_{(2)})\\
&=& \psi(f(h_{(2)})_{[0]}\ot h_{(1)}f(h_{(2)})_{[1]})\\
&=& f(h_{(2)})_{[0]}g(h_{(1)}f(h_{(2)})_{[1]}),
\end{eqnarray*}
as needed.
\end{proof}

The dual of the canonical map ${}^*\can:\ {}^*\Cc\to {}_A\Hom(A\ot_B A,A)$
is given by
$${}^*\can(\varphi)(b\ot_A a)=\varphi(\can(b\ot_A a))=\varphi(ba_{[0]}\ot a_{[1]},$$
and transports to
${}^*\can:\ \ul{\Hom}(H,A)\to {}_B\End(A)^{\rm op}$, given by
\begin{equation}\eqlabel{2.3.3}
{}^*\can(f)(a)={}^*\can(\varphi)(1\ot_B a)=a_{[0]}f(a_{[1]}).
\end{equation}

\begin{remarks}\relabel{2.4}
1) We have a projection $p:\ \Hom(H,A)\to \ul{\Hom}(H,A)$, $p(f)=\tilde{f}$,
given by 
$$\tilde{f}(h)=1_{[0]}f(1_{[1]}).$$
2) The unit element of $\ul{\Hom}(H,A)$ is not $\varepsilon$, but
$\tilde{\varepsilon}=\alpha(\varepsilon_\Cc)$. This can be verified directly
as follows:
\begin{eqnarray*}
&&\hspace*{-2cm}
(\tilde{\varepsilon}\# g)(h)=1_{[0]}\varepsilon(h_{(2)}1_{[2]})g(h_{(1)}1_{[1]})\\
&=& 1_{[0]}\varepsilon(h_{(2)}1_{(2)})g(h_{(1)}1_{(1)}1_{[1]})
= 1_{[0]}\varepsilon(h_{(2)}g(h_{(1)}1_{[1]})\\
&=& 1_{[0]}g(h1_{[1]})=g(h);\\
&&\hspace*{-2cm}
(f\# \tilde{\varepsilon})(h)=f(h_{(2)})_{[0]}\tilde{\varepsilon}
(h_{(1)}f(h_{(2)})_{[1]})\\
&=& f(h_{(2)})_{[0]}1_{[0']}{\varepsilon}(h_{(1)}f(h_{(2)})_{[1]}1_{[0']})=
f(h_{(2)})_{[0]}{\varepsilon} (h_{(1)}f(h_{(2)})_{[1]})\\
{\rm \equref{1.7.2}}~~~&=& 1_{[0]}\varepsilon(h_{(1)}1_{[1]})f(h_{(2)})\\
{\rm \equref{1.7.1}}~~~&=& 1_{[0]}f(h1_{[1]})=f(h).
\end{eqnarray*}
3) ${}^*\Cc$ is an $A$-bimodule, hence $\Hom(H,A)$ is also an $A$-bimodule.
The structure is given by the formula
\begin{equation}\eqlabel{2.4.1}
(afb)(h)=a_{[0]}f(ha_{[1]})b.
\end{equation}
\end{remarks}

Let us next compute the Morita context $(T,\ul{\Hom}(H,A),A,Q,\tau,\mu)$
from \prref{1.4}. We have a ring
morphism $j:\ A\rTo^{i} _{}{}^*\Cc\rTo^{\alpha} \ul{\Hom}(H,A)$, given by
$$j(a)(h)=1_{[0]}\varepsilon(h1_{[1]})a.$$
Take $\varphi\in Q=\{\varphi\in {}^*\Cc~|~c_{(1)}\varphi(c_{(2)})=
\varphi(c)\rho(1),~{\rm for~all~}c\in \Cc\}$, and let $\alpha(\varphi)=q$. Then
$$
c_{(1)}\varphi(c_{(2)})=(1_{[0]}\ot h_{(1)}1_{[1]})f(h_{(2)})=
f(h_{(2)})_{[0]}\ot h_{(1)}f(h_{(2)})_{[1]},$$
and
$$\varphi(c)\rho(1)=f(h)1_{[0]}\ot 1_{[1]},$$
hence $Q$ is the subset of $\ul{\Hom}(H,A)$ consisting of maps $f:\ H\to A$
satisfying
\begin{equation}\eqlabel{2.5.0}
f(h_{(2)})_{[0]}\ot h_{(1)}f(h_{(2)})_{[1]}=
f(h)1_{[0]}\ot 1_{[1]},
\end{equation}
for all $h\in H$.
$Q$ is an $(\ul{\Hom}(H,A),T)$-bimodule: $f\cdot q\cdot a= f\#q\# j(a)$,
for all $f\in \ul{\Hom}(H,A)$, $q\in Q$ and $a\in A$.\\
$A$ is a $(T,\ul{\Hom}(H,A))$-bimodule; the left $T$-action is given by
left multiplication, and the right $\ul{\Hom}(H,A))$-action is given by
$a\cdot q=a_{[0]}q(a_{[1]})$. The connecting maps
$$\tau:\ A\ot_{\ul{\Hom}(H,A)} Q\to T~~{\rm and}~~\mu:\ Q\ot_T A\to \ul{\Hom}(H,A)$$
are given by the formulas
$$\tau(a\ot q)=a_{[0]}q(a_{[1]});$$
$$\mu(q\ot a)=q\# j(a),~{\rm that~is~}\mu(q\ot a)(h)=q(h)a.$$
It follows from \prref{1.3} that $\tau$ is surjective if and only if there exists
$q\in Q$ such that $1_{[0]}q(1_{[1]})=1$. \thref{1.5} takes the following
form.

\begin{theorem}\thlabel{2.5}
Let $H$ be a finitely generated projective weak bialgebra, and $A$ a right
$H$-comodule algebra. Take a subring $B$ of $T=A^{{\rm co}H}$. Then,
with notation as above, the following
assertions are equivalent.
\begin{enumerate}
\item \begin{itemize}
\item ${\rm can}$ is an isomorphism;
\item $A$ is faithfully flat as a left $B$-module.
\end{itemize}
\item \begin{itemize}
\item ${}^*{\rm can}$ is an isomorphism;
\item $A$ is a left $B$-progenerator.
\end{itemize}
\item \begin{itemize}
\item $B=T$;
\item the Morita context $(B,\ul{\Hom}(H,A),A,Q,\tau,\mu)$ is strict.
\end{itemize}
\item \begin{itemize}
\item $B=T$;
\item $(F,G)$ is an equivalence of categories.
\end{itemize}
\end{enumerate}
\end{theorem}

\begin{proposition}\prlabel{2.6}
Let $H$ be a weak Hopf algebra. Then $H$ is a weak $H$-Galois extension of
$H^{{\rm co}H}=H^L$.
\end{proposition}

\begin{proof}
Let us first show that $H^{{\rm co}H}=H^L$. If $h\in H^{{\rm co}H}$, then
$\Delta(h)= h1_{(1)}\ot 1_{(2)}$, hence
$h=\varepsilon(h1_{(1)}) 1_{(2)}= \ol{\Pi}^R(h)\in H^L$.\\
Conversely, if $h\in H^L$, then $h=\ol{\Pi}^R(h)=\varepsilon(h1_{(1)}) 1_{(2)}$,
so 
\begin{eqnarray*}
&&\hspace*{-2cm}
\Delta(h)=\varepsilon(h1_{(1)})1_{(2)}\ot 1_{(3)}\\
&=& \varepsilon(h1_{(1)})1_{(2)}1_{(1')}\ot 1_{(2')}= h1_{(1)}\ot 1_{(2)},
\end{eqnarray*}
hence $h\in H^{{\rm co}H}$.\\
Let us next show that $\can$ is invertible, with inverse
$$\can^{-1}(1_{(1)}\ot  h1_{(2)})=1_{(1)}S(h_{(1)}1_{(2)}\ot_{H^L} h_{(2)}1_{(3)}.$$
\begin{eqnarray*}
&&\hspace*{-2cm}
\can^{-1}(\can (h\ot k))= \can^{-1}(hk_{(1)}\ot k_{(2)})= hk_{(1)}S(k_{(2)})\ot_{H^L}  k_{(3)}\\
&=& h\Pi^L(k_{(1)})\ot_{H^L}  k_{(2)}
= h\ot_{H^L}\Pi^L(k_{(1)})  k_{(2)}\\
&=& h\ot_{H^L}\varepsilon(1_{(1)}k_{(1)})1_{(2)}k_{(2)}=
h\ot_{H^L}(k);\\
\end{eqnarray*}
\begin{eqnarray*}
&&\hspace*{-2cm}
\can(\can^{-1}(1_{(1)}\ot  h1_{(2)}))=
\can(1_{(1)}S(h_{(1)}1_{(2)}\ot_{H^L} h_{(2)}1_{(3)})\\
&=& 1_{(1)}S(h_{(1)}1_{(2)})h_{(2)}1_{(3)}\ot h_{(3)}1_{(4)}\\
&=& 1_{(1)}1_{(1')}\varepsilon(h_{(1)}1_{(2)}1_{(2')})\ot h_{(2)}1_{(3)}\\
&=& 1_{(1)}1_{(1')}\varepsilon(h_{(1)}1_{(3)})\varepsilon(1_{(2)}1_{(2')})\ot h_{(2)}1_{(4)}\\
&=& 1_{(1)}1_{(1')}\varepsilon(1_{(2)}1_{(2')})\ot \varepsilon(h_{(1)}1_{(3)})h_{(2)}1_{(4)}\\
&=&
1_{(1)}1_{(1')}\varepsilon(1_{(2)}1_{(2')})\ot \varepsilon h1_{(3)})
= 1_{(1)}\ot  h1_{(2)}.
\end{eqnarray*}
\end{proof}

\section{Groupoid gradings}\selabel{3}
Recall that a groupoid $G$ is a category in which every morphism is an isomorphism.
In this Section, we consider finite groupoids, i.e. groupoids with a finite number of
objects. The set of objects of $G$ will be denoted by $G_0$, and the set of morphisms
by $G_1$. The identity morphism on $x\in G_0$ will also be denoted by $x$.
For $\sigma:\ x\to y$ in $G_1$, we write
$$s(\sigma)=x~~{\rm and}~~t(\sigma)=y,$$
respectively for the source and the target of $\sigma$. For every
$x\in G$, $G_x=\{\sigma\in G~|~s(\sigma)=t(\sigma)=x\}$ is a group.\\
Let $G$ be a groupoid, and $k$ a commutative ring. The groupoid algebra is the
direct product
$$kG=\bigoplus_{\sigma\in G_1}ku_{\sigma},$$
with multiplication defined by the formula
$$u_\sigma u_\tau=\begin{cases}
u_{\sigma\tau}&{\rm if }~t(\tau)=s(\sigma);\\
0& {\rm if } ~t(\tau)\neq s(\sigma).
\end{cases}$$
The unit element is $1=\sum_{x\in G_0} u_x$. $kG$ is a weak Hopf algebra, with
comultiplication, counit and antipode given by the formulas
$$\Delta(u_\sigma)= u_\sigma \ot u_\sigma,~\varepsilon(u_\sigma)=1~{\rm and}~
S(u_\sigma)=u_{\sigma^{-1}}.$$
Using the formula
$$\Delta(1)=\sum_{x\in G_0} u_x\ot u_x.$$
We compute that $\Pi^L:\ kG\to kG$ is given by the formula
$$\Pi^L(u_\sigma)=\sum_{x\in G_0}\varepsilon(u_x u_\sigma)=u_{t(\sigma)},$$
hence
$$(kG)^L= \bigoplus_{x\in G_0} ku_x.$$
Let $k$ be a commutative ring. A $G$-graded $k$-algebra is a $k$-algebra $A$
together with a direct sum decomposition
$$A=\bigoplus_{\sigma\in G_1} A_\sigma,$$
such that
\begin{equation}\eqlabel{3.1.1}
A_\sigma A_\tau
\begin{cases}
\subset A_{\sigma\tau}&{\rm if}~~t(\tau)=s(\sigma);\\
=0&{\rm if}~~t(\tau)\neq s(\sigma).
\end{cases}
\end{equation}
and
\begin{equation}\eqlabel{3.1.2}
1_A\in \bigoplus_{x\in G_0}A_x.
\end{equation}

\begin{proposition}\prlabel{3.1}
Let $G$ be a finite groupoid, and $k$ a commutative ring. We have an isomorphism
between the categories of $kG$-comodule algebras and $G$-graded $k$-algebras.
\end{proposition}

\begin{proof}
Let $(A,\rho)$ be a $kG$-comodule algebra, and define $\rho:\ A\to A\ot kG$
by
$$A_\sigma=\{a\in A~|~\rho(a)=a\ot u_\sigma\}.$$
From the fact that $A\ot kG$ is a free left $A$-module with basis
$\{u_\sigma~|~\sigma\in G_1\}$, it follows that $A_\sigma\cap A_\tau=\{0\}$
if $\sigma\neq \tau$.\\
For $a\in A$, we can write
$$\rho(a)=\sum_{\sigma\in G_1} a_\sigma\ot u_\sigma.$$
From the coassociativity of $\rho$, it follows that
$$\sum_{\sigma\in G_1} \rho(a_\sigma)\ot u_\sigma=
\sum_{\sigma\in G_1} a_\sigma\ot u_\sigma\ot u_\sigma,$$
hence $\rho(a_\sigma)a_\sigma\ot u_\sigma$, so $a_\sigma\in A_\sigma$ and
$$a=\sum_{\sigma\in G_1} a_\sigma\varepsilon( u_\sigma)=
\sum_{\sigma\in G_1} a_\sigma\in _{\sigma\in G_1} A_\sigma.$$
If $a\in A_\sigma$ and $b\in A_\tau$, then $\rho(ab)=ab\ot u_\sigma u_\tau$,
and \equref{3.1.1} follows.\\
\equref{1.6.9} tells us that $\rho(1_A)\in A\ot (kG)^L=\bigoplus_{x\in G_0}
A\ot u_x$, hence we can write
$$\rho(1_A)=\sum_{x\in G_0} 1_x\ot u_x,$$
and
$$1_A=\sum_{x\in G_0} 1_x\in \bigoplus_{x\in G_0}A_x.$$
Conversely, let $A$ be a $G$-graded algebra. For $a\in A$, let
$a_\sigma$ be the projection of $A$ on $A_\sigma$. Then define $\rho(a)=
\sum_{\sigma\in G_1} a_\sigma\ot u_\sigma$. \equref{1.6.9} then follows
immediately from \equref{3.1.2}.
\end{proof}

\begin{proposition}\prlabel{3.2}
Let $A$ be a $G$-graded algebra, and take $x\in G_0$. Then
$\bigoplus_{\sigma\in G_x} A_\sigma$ is a $G_x$-graded algebra, with
unit $1_x$. In particular, $A_x$ is a $k$-algebra with unit $1_x$.
\end{proposition}

\begin{proof}
If $a\in A_\sigma$, with $t(\sigma)=x$, then
$$a\ot \sigma=\rho(a)=\rho(1_Aa)=\sum_{y\in G_0} 1_ya\ot u_yu_\sigma=
1_xa\ot \sigma,$$
hence $a=1_x a$. In a similar way, we can prove that $a=a1_x$ if $s(a)=x$.
The result then follows easily.
\end{proof}

Let $A$ be a $G$-graded algebra. A $G$-graded right $A$-module is a
right $A$-module together with a direct sum
decomposition
$$M=\bigoplus_{\sigma\in G_1}M_\sigma$$
such that
$$M_\sigma A_\tau
\begin{cases}
\subset M_{\sigma\tau}&{\rm if}~~t(\tau)=s(\sigma);\\
=0&{\rm if}~~t(\tau)\neq s(\sigma).
\end{cases}$$
A right $A$-linear map $f:\ M\to N$ between two $G$-graded right $A$-modules
is called graded if $f(M_\sigma)\subset N_\sigma$, for all $\sigma\in G_1$.
The category of $G$-graded right $A$-modules and graded $A$-linear maps
is denoted by $\Mm_A^G$.
The proof of the next result is then similar to the proof of \prref{3.1}.

\begin{proposition}\prlabel{3.3}
Let $G$ be a groupoid, and $A$ a $G$-graded algebra. Then the categories
$\Mm_A^G$ and $\Mm_A^{kg}$ are isomorphic.
\end{proposition}

If $m\in M_\sigma$, with $s(\sigma)=x$, then
$$\rho(m)=m\ot \sigma= \rho(m1)=\sum_{y\in G_0}m1_y\ot \sigma x=
m1_x\ot \sigma,$$
hence $m1_x=m$.\\

Let $M\in \Mm_A^G\cong \Mm_A^{kg}$.  Then $m\in M^{{\rm co}kG}$ if and only if
$$\rho(m)=\sum_{\sigma\in G_1} m_\sigma\ot\sigma=\sum_{x\in G_0}m1_x\ot x,$$
if and only if $m\in \bigoplus_{x\in G_0} M_x$. We conclude that
\begin{equation}\eqlabel{3.4.1}
M^{{\rm co}kG}= \bigoplus_{x\in G_0} M_x.
\end{equation}
In particular,
$$T=A^{{\rm co}kG}= \bigoplus_{x\in G_0} A_x.$$
If $N\in \Mm_T$, then $N=\bigoplus_{x\in G} N_x$, with
$N_x=N1_x\in \Mm_{A_x}$.
We have a pair of adjoint functors
$$F=-\ot_T A:\ \Mm_T\to \Mm_A^G,~~G=(-)^{{\rm co}kG}:\ \Mm_A^G\to \Mm_T,$$
with unit and counit given by the formulas
$$\nu_N:\ N\to (N\ot_T A)^{{\rm co}kG},~~\nu_N(n)=\sum_{x\in G_0} n_x\ot_T 1_x$$
$$\zeta_M:\ M^{{\rm co}kG}\ot_T A\to M,~~\zeta_M(m\ot_T a)=ma.$$
A $G$-graded $k$-algebra is called strongly graded if
$$A_\sigma A_\tau=A_{\sigma\tau}~~{\rm if}~~t(\tau)=s(\sigma).$$

\begin{proposition}\prlabel{3.4}
Let $A$ be a $G$-graded $k$-algebra. With notation as above,
$\nu_N:\ N\to (N\ot_T A)^{{\rm co}kG}$ is bijective, for every $N\in \Mm_T$.
\end{proposition}

\begin{proof}
 First observe that
$$(N\ot_T A)_\sigma= N\ot_T A_\sigma= N_{t(\sigma)}\ot_{T} A_\sigma,$$
hence
$$(N\ot_T A)^{{\rm co}kG}=\bigoplus_{x\in G_0} (N\ot_T A)_x=
\bigoplus_{x\in G_0} N_x\ot_{T} A_x.$$
The map
$$f:\ (N\ot_T A)^{{\rm co}kG}\to N,~~f(n_x\ot a_x)=n_xa_x$$
is the inverse of $\nu_N$: it is obvious that
$f\circ \nu_N=N$. We also have that
$$(\nu_N\circ f)(n_x\ot_{T} a_x)=n_xa_x\ot_{T} 1_x=n_x\ot_{T} a_x.$$
\end{proof}

Let $A$ be an algebra graded by a group $G$. It is a classical result from graded ring theory
(see e.g. \cite{NVO}) that $A$ is strongly graded if and only if 
the categories of $G$-graded $A$-modules and $A_1$-modules are equivalent.
This is also equivalent to $A$ being a $kG$-Galois extension of $A_e$.
We now present the groupoid version of this result.

\begin{theorem}\thlabel{3.5}
Let $G$ be a groupoid, and $A$ a $G$-graded $k$-algebra. Then the following
assertions are equivalent.
\begin{enumerate}
\item $A$ is strongly graded;
\item $(F,G)$ is a category equivalence;
\item $(A\ot kG, \rho(1)=\sum_{x\in G_0} e_x\ot x)$ is a Galois coring;
\item $\can:\ A\ot_T A\to A\ot kG,~~\can(a\ot b)=\sum_{\sigma\in G_1}ab_\sigma
\ot \sigma$ is surjective.
\end{enumerate}
\end{theorem}

\begin{proof}
$\ul{1)\Rightarrow 2)}$. In view of \prref{3.4}, we only have to show
that $\zeta_M$ is bijective, for every graded $A$-module $M$.\\
Take $m\in M_\sigma$, with $s(\sigma)=x$, $t(\sigma)=y$. Then
$\sigma\sigma^{-1}=x$, and there exist $a'_i\in A_{\sigma^{-1}}$,
$a_i\in A_\sigma$ such that $\sum_i a'_ia_i=1_x$. We have that
$ma'_i\in M_\sigma A_{\sigma^{-1}}\subset M_y\subset M^{{\rm co}kG}$,and
$$\zeta_M(\sum_i ma'_i\ot a_i)=\sum_i ma'_ia_i=m,$$
and it follows that $\zeta_M$ is surjective.\\
Take 
$$m_j=\sum_{x\in G_0} m_{j,x}\in \bigoplus_{x\in G_0} M_x~~{\rm and}~~
c_j\in A.$$
Assume
$$\zeta_M(\sum_j m_j\ot_T c_j)=\sum_j m_jc_j=0,$$
and take $\sigma\in G_1$ with $s(\sigma)=x$ and $t(\sigma)=y$. Then
$$0=(\sum_j m_jc_j)_\sigma=\sum_j m_{j}c_{j,\sigma},$$
hence
$$\sum_j m_{j}\ot_T c_{j,\sigma}=\sum_{i,j} m_{j}\ot_T c_{j,\sigma}a'_ia_i=
=\sum_{i,j} m_{j} c_{j,\sigma}a'_i\ot_T a_i=0,$$
and
$$\sum_j m_j\ot_T c_j=\sum_{\sigma\in G_1}\sum_j m_{j}\ot_T c_{j,\sigma}=0,$$
and it follows that $\zeta_M$ is injective.\\
$\ul{2)\Rightarrow 3)}$ follows from the observation made before \thref{1.2}
(see \cite[Prop. 3.1]{Caenepeel03}).\\
$\ul{3)\Rightarrow 4)}$ is trivial.\\
$\ul{4)\Rightarrow 1)}$. Take $\sigma,\tau\in G_1$ such that
$s(\sigma)=t(\tau)$, and $c\in A_{\sigma\tau}$. Since $\can$ is surjective,
there exist homgeneous $a_i, b_i\in A$ such that
$$\can(\sum_{i} a_i\ot b_i)=\sum_{i,j}a_i b_j\ot {\rm deg}(b_j)
=c\ot \tau.$$
On the left hand side, we can delete all the terms for which ${\rm deg}(b_j)
\neq \tau$. So we find
$$\sum_{i,j}a_i b_j\ot \tau
=c\ot \tau,$$
and
$$\sum_{i,j}a_i b_j=c.$$
On the left hand side, we can now delete all terms of degree different from
$\sigma\tau$, since the degree of the right hand side is $\sigma\tau$.
This means that we can delete all terms for which $\deg(a_i)\neq \sigma$.
So we find that $\sum_{i,j}a_i b_j=c$, with $\deg(a_i)=\sigma$, and
$\deg(b_i)=\tau$, and $A_\sigma A_\tau=A_{\sigma\tau}$.
\end{proof}

\section{Groupoid actions}\selabel{4}
Let $G$ be a groupoid, as in \seref{3}, with the additional assumption that
$G_1$ is finite. Then $kG$ is free of finite rank as a $k$-module, hence
$Gk=(kG)^*$ is also a weak Hopf algebra. As a $k$-module,
$$Gk=\bigoplus_{\sigma\in G_1} kv_\sigma,$$
with
$$\lan v_\sigma,\tau\ran=\delta_{\sigma,\tau}.$$
The algebra structure is given by the formulas
$$v_\sigma v_\tau=\delta_{\sigma,\tau}v_\sigma~~;~~1=\sum_{\sigma\in G_1}v_\sigma,$$
and the coalgebra structure is
$$\Delta(v_\sigma)=\sum_{\tau\rho=\sigma} v_\tau\ot v_rho=
\sum_{t(\tau)=t(\sigma)} v_\tau\ot v_{\tau^{-1}\sigma},$$
$$\varepsilon(\sum_{\sigma\in G_1}a_\sigma v_\sigma)=
\sum_{x\in G_0}a_x v_x.$$
The antipode is given by $S(v_{\sigma})=v_{\sigma^{-1}}$. Observe that
$$\Delta(1)=1_{(1)}\ot 1_{(2)}=\sum_{t(\rho)=s(\sigma)} v_\tau\ot v_\sigma.$$
Let $A$ be a right $Gk$-comodule algebra. Then $A$ is a left $kG$-module algebra,
with left $kG$-action given by
$$\sigma\cdot a=\lan a_{[1]},\sigma\ran a_{[0]}.$$
The $Gk$-coaction $\rho$ can be recovered from the action, using the formula
\begin{equation}\eqlabel{4.1.1}
\rho(a)=\sum_{\sigma\in G_1}\sigma\cdot a\ot v_\sigma.
\end{equation}
In \cite[Prop. 4.15]{CDG}, equivalent definitions of an $H$-module algebra
are given. If we apply them in the case where $H=kG$, we find that a
$k$-algebra with a left $kG$-module structure is a left $kG$-module algebra
if 
\begin{equation}\eqlabel{4.1.2}
\sigma\cdot (ab)=(\sigma\cdot a)(\sigma\cdot b).
\end{equation}
 for all $\sigma\in G_1$
and $a,b\in A$, and the following equivalent conditions are satisfied.
\begin{eqnarray}
&&\sigma\cdot 1_A=t(\sigma)\cdot 1_A\eqlabel{4.1.3}\\
&&a(\sigma\cdot 1_A)= t(\sigma)\cdot a\eqlabel{4.1.4}\\
&&(\sigma\cdot 1_A)a= t(\sigma)\cdot a\eqlabel{4.1.5}
\end{eqnarray}
We will call $A$ a left $G$-module algebra.

\begin{proposition}\prlabel{4.1}
Let $A$ be a left $G$-module algebra. Then $\{x\cdot 1_A~|~x\in G_0\}$
is a set of central orthogonal idempotents in $A$. Hence we can write
$$A=\bigoplus_{x\in G_0} A_x,$$
where $A_\sigma=(\sigma\cdot 1_A)A=(t(\sigma)\cdot 1_A)=A_{t(\sigma)}$.
\end{proposition}

\begin{proof}
It follows from \equref{4.1.1}, with $a=b=1_A$ and $\sigma=x$, that 
$x\cdot 1_A$ is idempotent. From \equref{4.1.4} and \equref{4.1.5},
it follows that $x\cdot 1_A$ is central.
If $x\neq y\in G_0$, then we find, using \equref{4.1.4}, that
$(x\cdot 1_A)(y\cdot 1_A)=y\cdot(x\cdot 1_A)=yx\cdot 1_A=0$.
Finally
$$\sum_{x\in G_0} x\cdot 1_A=(\sum_{x\in G_0}x)\cdot 1_A=1_{kG}\cdot 1_A=1_A.$$
Observe also that $\sigma\cdot a=(\sigma\cdot a)(\sigma\cdot 1_A)\in A_{\sigma}$.
\end{proof}

We have
$$A\ot Gk=\bigoplus_{\sigma\in G_1} Av_\sigma.$$
Let us compute the projection
$$g:\ \bigoplus_{\sigma\in G_1} Av_\sigma\to \bigoplus_{\sigma\in G_1} Av_\sigma$$
introduced in \seref{2}.
$$
g(av_\sigma)= \sum_{\tau\in G_1} a(\tau\cdot 1_A)\ot v_\sigma v_\tau
= a(\sigma\cdot 1_A)\ot v_\sigma.$$
We have seen in \leref{2.1} that we have an $A$-coring
$$\Cc=\im(g)=\bigoplus_{\sigma\in G_1} A_\sigma v_\sigma.$$
The right $A$-module structure is given by the formula
$$((\sigma\cdot 1_A)v_\sigma)a=(\sigma\cdot a)v_\sigma.$$
Now assume that $M\in \Mm^\Cc$, or, equivalently, $M\in \Mm^{Gk}_A$
(see \prref{2.2}). Then $M$ is a right $A$-module, and also a right $Gk$-comodule,
and a fortiori a left $kG$-module. Condition \equref{2.2.1} is then equivalent to
$$\sum_{\sigma\in G_1} \sigma\cdot(ma)\ot v_\sigma=
\sum_{\sigma,\tau\in G_1} (\sigma\cdot m)(\tau\cdot a)\ot v_\sigma v_\tau=
\sum_{\sigma\in G_1}(\sigma\cdot m)(\sigma\cdot a)\ot v_\sigma,$$
or
\begin{equation}\eqlabel{4.2.1}
\sigma\cdot(ma)=(\sigma\cdot m)(\sigma\cdot a),
\end{equation}
for all $\sigma\in G_1$, $m\in M$ and $a\in A$. We will also say that $G$ acts
as a groupoid of right $A$-semilinear automorphisms on $M$. Now
$m\in M^G=M^{{\rm co}\Cc}$ if and only if 
$$\rho(m)=\sum_{\sigma\in G_1} \sigma\cdot m\ot v_\sigma$$
equals
$$m1_{[0]}\ot 1_{[1]}=\sum_{\sigma\in G_1}
m(\sigma\cdot 1_A)\ot v_\sigma,$$
hence
$$M^G=\{m\in M~|~\sigma\cdot m= m(\sigma\cdot 1_A),~{\rm for~all~}\sigma\in G_1\}.$$
In particular,
$$T=A^G=\{a\in A~|~\sigma\cdot a=t(\sigma)\cdot a,~{\rm for~all~}\sigma\in G_1\}.$$
Let $B\to T$ be a ring morphism. Let us compute
$$\can:\ A\ot_B A\to \bigoplus_{\sigma\in G_1} A_\sigma v_\sigma.$$
\begin{eqnarray*}
\can(a\ot b)&=&
\sum_{\sigma\in G_1} (a(\sigma \cdot 1_A)\ot v_\sigma)b\\
&=& \sum_{\sigma\in G_1} a(\sigma \cdot 1_A)(\sigma \cdot b)\ot v_\sigma\\
&=& \sum_{\sigma\in G_1} a(\sigma \cdot b)\ot v_\sigma.
\end{eqnarray*}

Our next goal is to compute 
$${}^*\Cc=\ul{\Hom}(Gk,A)=\{f:\ Gk\to A~|~1_{[0]}f(v_\tau 1_{[1]})=f(v_\tau),~
{\rm for~all}~\tau\in G_1\}.$$

\begin{proposition}\prlabel{4.3}
Let $G$ be finite groupoid, and $A$ a $G$-module algebra. Then
$$\ul{\Hom}(Gk,A)=\bigoplus_{\sigma\in G_1} u_\sigma A_\sigma,$$
as a right $A$-modules. The left $A$-module structure on $\ul{\Hom}(Gk,A)$
is given by the formula
\begin{equation}\eqlabel{4.3.1}
a u_\sigma(\sigma\cdot 1_A)=u_\sigma(\sigma\cdot a),
\end{equation}
and the multiplication is given by
\begin{equation}\eqlabel{4.3.2}
U_\sigma\# U_\tau=U_{\tau\sigma},
\end{equation}
where we denoted $U_\sigma=u_\sigma(\sigma\cdot 1_A)$. The unit element of
$\ul{\Hom}(Gk,A)$ is 
$$\sum_{x\in G_0} U_x.$$
\end{proposition}

\begin{proof}
Observe that $\Hom(Gk,A)\cong (Gk)^*\ot A=kG\ot A$, so that every $f:\ Gk\to A$
can be written as
$$f=\sum_{\sigma\in G_1}u_\sigma a_\sigma,$$
with $a_\sigma=f(v_\sigma)$. $u_\sigma$ is then the projection onto the
component $v_\sigma$ of $Gk$. We find that $f\in \ul{\Hom}(Gk,A)$ if and only
if $f(v_\tau)=a_\tau$ equals
$$ 1_{[0]}f(v_\tau 1_{[1]})=
\sum_{\sigma\in G_1} (\sigma\cdot 1_A)f(v_\tau v_\sigma)=(\tau\cdot 1_A)a_\tau,$$
or, equivalently, $a_\tau\in A_\tau$, for every $\tau\in G_1$.
Write $U_\sigma=u_\sigma(\sigma\cdot 1_A)$, then we have
$$\ul{\Hom}(Gk,A)=\bigoplus_{\sigma\in G_1} U_\sigma A=\bigoplus_{\sigma\in G_1} U_\sigma A_\sigma.$$
Observe that $U_\sigma(v_\tau)=(\sigma\cdot 1_A)\delta_{\sigma,\tau}$.\\
The left $A$-action on $\ul{\Hom}(Gk,A)$ is computed using \equref{2.4.1}. We compute
\begin{eqnarray*}
(aU_\sigma)(v_\tau)&=&
a_{[0]}U_\sigma(v_\tau a_{[1]})=\sum_{\lambda\in G_1}
(\lambda\cdot a)U_\sigma(v_\tau v_\lambda)\\
&=& (\tau\cdot a)U_\sigma(v_\tau)=(\tau\cdot a)(\sigma\cdot 1_A)\delta_{\sigma,\tau},
\end{eqnarray*}
hence $(aU_\sigma)(v_\tau)=0$ if $\sigma\neq \tau$. If $\sigma=\tau$, then
$$(aU_\sigma)(v_\sigma)=(\sigma\cdot a)(\sigma\cdot 1_A)
=(\sigma\cdot a)= (U_\sigma(\sigma\cdot a))(v_\sigma),$$
and we conclude that
$$aU_\sigma=U_\sigma(\sigma\cdot a).$$
Using \equref{2.3.2}, we compute
\begin{eqnarray*}
&&\hspace*{-2cm}
(U_\sigma\# U_\tau)(v_\lambda)=
\sum_{\mu\nu=\lambda}U_\sigma(v_\nu)_{[0]}U_\tau(v_\mu f(v_\nu)_{[1]})\\
&=& \sum_{\rho\in G_1}\sum_{\mu\nu=\lambda}
\rho\cdot((\sigma\cdot 1_A)\delta_{\sigma,\nu})U_\tau(v_\mu v_\rho)\\
&=& \sum_{\mu\nu=\lambda} \mu((\sigma\cdot 1_A)\delta_{\sigma,\nu})(\tau\cdot 1_A)
\delta_{\tau,\mu}.
\end{eqnarray*}
If $\tau\sigma\neq \lambda$, then $(U_\sigma\# U_\tau)(v_\lambda)=0$.\\
If $\tau\sigma= \lambda$, then
\begin{eqnarray*}
&&\hspace*{-2cm}
(U_\sigma\# U_\tau)(v_\lambda)=(\tau\cdot(\sigma\cdot 1_A))(\tau\cdot 1_A)
=(\tau\sigma\cdot 1_A)(\tau\cdot 1_A)\\
&=& (t(\tau\sigma)\cdot 1_A)(t(\tau)\cdot 1_A)
= (t(\tau\sigma)\cdot 1_A)(t(\tau\sigma)\cdot 1_A)\\
&=& t(\tau\sigma)\cdot 1_A=(\tau\sigma)\cdot 1_A=U_{\tau\sigma}(v_\lambda),
\end{eqnarray*}
and we conclude that $U_\sigma\# U_\tau=U_{\tau\sigma}$.
\end{proof}

Recall that a ring morphism $A\to R$ is called {\sl Frobenius}
if there exists an $A$-bimodule map $\ol{\nu}:\ R\to A$ and
$e=e^1\ot_A e^2\in R\ot_AR$ (summation implicitly understood) such that
\begin{equation}\eqlabel{4.4.1}
re^1\ot_A e^2=e^1\ot_A e^2r
\end{equation}
for all $r\in R$, and
\begin{equation}\eqlabel{4.4.2}
\ol{\nu}(e^1)e^2=e^1\ol{\nu}(e^2)=1.
\end{equation}
This is equivalent to the restrictions of scalars $\Mm_R\to
\Mm_A$ being a Frobenius functor, which means that its left and right
adjoints are isomorphic (see \cite[Sec. 3.1 and 3.2]{CMZ}).
$(e,\ol{\nu})$ is then called a Frobenius system.

\begin{proposition}\prlabel{4.4}
Let $G$ be finite groupoid, and $A$ a $G$-module algebra. Then
the ring morphism $A\to \ul{\Hom}(Gk,A)={}^*\Cc$ is Frobenius.
\end{proposition}

\begin{proof}
The Frobenius system is the following:
$$e=\sum_{\sigma\in G_1} U_{\sigma^{-1}}\ot U_\sigma;$$
$$\ol{\nu}(\sum_{\sigma\in G_1})u_\sigma (\sigma\cdot 1_A)a_\sigma)=
\sum_{x\in G_0}(x\cdot 1_A)a_x.$$
For all $a\in A$, we have that
\begin{eqnarray*}
ae&=& \sum_{\sigma\in G_1} U_{\sigma^{-1}}(\sigma^{-1}\cdot a)\ot U_\sigma
= \sum_{\sigma\in G_1} U_{\sigma^{-1}}\ot U_\sigma
((\sigma\sigma^{-1})\cdot a)\\
&=& \sum_{\sigma\in G_1} U_{\sigma^{-1}}\ot U_\sigma
(t(\sigma)\cdot a)
= \sum_{\sigma\in G_1} U_{\sigma^{-1}}\ot U_\sigma
(\sigma\cdot 1_A)a=ea;
\end{eqnarray*}
and
\begin{eqnarray*}
&&\hspace*{-2cm}
\ol{\nu}(a\sum_{\sigma\in G_1})u_\sigma (\sigma\cdot 1_A)a_\sigma)
=\ol{\nu}(\sum_{\sigma\in G_1})u_\sigma (\sigma\cdot a)a_\sigma)\\
&=& \sum_{x\in G_0} (x\cdot a)a_x= \sum_{x\in G_0} a(x\cdot 1_A)a_x\\
&=& a\ol{\nu}(\sum_{\sigma\in G_1})u_\sigma (\sigma\cdot 1_A)a_\sigma),
\end{eqnarray*}
and it follows that $\ol{\nu}$ is left $A$-linear. It is clear that $\ol{\nu}$
is right $A$-linear. Finally,
\begin{eqnarray*}
\ol{\nu}(e^1)e^2&=& \sum_{x\in G_0} (x\cdot 1_A)u_x(x\cdot 1_A)\\
&=& \sum_{x\in G_0} (x\cdot 1_A)u_x(x^2\cdot 1_A)(x\cdot 1_A)=\sum_{x\in G_0} U_x=
1.
\end{eqnarray*}
In a similar way, we show that $e^1\ol{\nu}(e^2)=1$.
\end{proof}

Using \equref{2.3.3}, we compute the map ${}^*\can:\ \ul{\Hom}(H,A)\to
{}_B\End(A)$:
$${}^*\can(U_\sigma a_\sigma)(b)=(\sigma\cdot b)a_\sigma.$$

Let us now compute the module $Q$ (see \equref{2.5.0}).
From \cite[Theorem 2.7]{CVW}, it follows that $Q$ and $A$ are isomorphic
as abelian groups. This will follow from our computations.

\begin{proposition}\prlabel{4.5}
Let $G$ be finite groupoid, and $A$ a $G$-module algebra. Then
$$Q=\{\sum_{\sigma\in G_1} u_{\sigma} (\sigma\cdot a)~|~a\in A\}.$$
\end{proposition}

\begin{proof}
Take $f=\sum_{\sigma\in G_1} U_\sigma a_\sigma\in \Hom(H,A)$. Then $f\in Q$
if and only if \equref{2.5.0} is satisfied.
Take $h=v_\tau$ in \equref{2.5.0}.
The right hand side of \equref{2.5.0} is
$$\sum_{\mu\in G_1} a_\tau (\mu\cdot 1_A)\ot v_\mu=
\sum_{\mu\in G_1} (\mu\cdot 1_A)a_\tau\ot v_\mu.$$
The left hand side of \equref{2.5.0} is
$$\sum_{\sigma\in G_1}\sum_{\mu\nu=\tau}\sum_{\rho\in G_1}
\rho\cdot (U_\sigma a_\sigma)(v_\nu)\ot v_\mu v_\rho=
\sum_{\mu\nu=\tau} \mu\cdot a_\nu\ot v_\mu.$$
It follows that $f\in Q$ if and only if
\begin{equation}\eqlabel{4.5.1}
\sum_{\mu\in G_1} (\mu\cdot 1_A)a_\tau\ot v_\mu=
\sum_{\mu\nu=\tau} \mu\cdot a_\nu\ot v_\mu,
\end{equation}
for all $\tau\in G_1$. Assume that $f\in Q$, and
take the $v_\tau$ component of \equref{4.5.1}. Then $\nu=s(\tau)$ and
\begin{equation}\eqlabel{4.5.2}
(\tau\cdot 1_A)a_\tau=a\tau=\tau\cdot a_{s(\tau)}.
\end{equation}
This means that $f$ is completely determined by 
$$a=\sum_{x\in G_0} a_x\in A=\bigoplus_{x\in G_0} A_x.$$
Conversely, take $a\in A$, let $a_x=a(x\cdot 1_A)\in A_x$, and define
$a_\tau$ using \equref{4.5.2}. Then
$$(\tau\cdot 1_A)a_\tau=(\tau\cdot 1_A)(\tau\cdot a_{s(\tau)})=
\tau\cdot(1_A s(\tau))a_\tau,$$
so $a_\tau\in A_\tau$, as needed. Also observe that
$$a_\tau=\tau\cdot a(s(\tau)\cdot 1_A))=
(\tau\cdot a)((\tau s(\tau))\cdot 1_A)=(\tau\cdot a)(\tau \cdot 1_A)=
\tau\cdot a.$$
We then claim that
$$f=\sum_{\sigma\in G_1} U_\sigma a_\sigma\in Q.$$
Take $\sigma,\tau\in G_1$ with $t(\sigma)\neq t(\tau)$. Then for all $a\in A$,
we have
$$(\sigma\cdot 1_A)(\tau\cdot 1_A)a=(\sigma\cdot 1_A)(t(\tau)\cdot a)
= (t(\sigma)t(\tau))\cdot a=0.$$
It follows that the left hand side of \equref{4.5.1} amounts to
$$\sum_{t(\mu)=t(\tau)} (\mu\cdot 1_A)a_\tau\ot v_\mu=
\sum_{t(\mu)=t(\tau)} a_\tau\ot v_\mu.$$
The right hand side of \equref{4.5.1} is
$$\sum_{t(\mu)=t(\tau)} \mu\cdot a_{\mu^{-1}\tau}\ot v_\mu.$$
If $t(\mu)=t(\tau)$, then
$$\mu\cdot a_{\mu^{-1}\tau}=\mu\cdot ((\mu^{-1}\tau)\cdot a_{s(\mu^{-1}\tau)}=
\tau\cdot a_{s(\tau)}= a_\tau,$$
and \equref{4.5.1} follows. Hence $f\in Q$.
\end{proof}

We have seen in \seref{2} that $Q$ is a $(\ul{\Hom}(Gk,A),T)$-bimodule.
Using \prref{4.5}, we can transport this bimodule structure to $A$.
We find the following bimodule structure on $A$:
$$(\sum_{\tau\in G_1}U_\tau b_\tau)\cdot a\cdot u=
\sum_{\tau\in G_1} \tau^{-1}\cdot(b_\tau a)u.$$
$A$ is also a $(T, \ul{\Hom}(Gk,A))$-bimodule. The structure is
$$u\bullet a\bullet (\sum_{\tau\in G_1}U_\tau b_\tau)=
\sum_{\tau\in G_1}u(\tau\cdot a)b_\tau.$$
We have a Morita context $(T,\ul{\Hom}(Gk,A),A,Q,\tau,\mu)$, and, a fortiori,
a Morita context $(T,\ul{\Hom}(Gk,A),A,A,\tau,\mu)$. The connecting maps
$$\tau:\ A\ot_{\ul{\Hom}(Gk,A)} A\to T~~{\rm and}~~\mu:\ A\ot_T A\to \ul{\Hom}(Gk,A)$$
are given by the formulas
$$\tau(b\ot a)=\sum_{\sigma\in G_1} \sigma\cdot (ba);$$
$$\mu(a\ot b)=\sum_{\sigma\in G_1}U_\sigma (\sigma\cdot a)b.$$
It follows from \prref{1.3} that $\tau$ is surjective if and only if there exists
$a\in A$ such that $\sum_{\sigma\in G_1} \sigma\cdot a=1$. \thref{2.5} takes the following form.

\begin{theorem}\thlabel{4.5}
Let $G$ be a finite groupoid, and $A$ a $G$-module algebra.
Take a subring $B$ of $T=A^{G}$. Then,
with notation as above, the following
assertions are equivalent.
\begin{enumerate}
\item \begin{itemize}
\item ${\rm can}$ is an isomorphism;
\item $A$ is faithfully flat as a left $B$-module.
\end{itemize}
\item \begin{itemize}
\item ${}^*{\rm can}$ is an isomorphism;
\item $A$ is a left $B$-progenerator.
\end{itemize}
\item \begin{itemize}
\item $B=T$;
\item the Morita context $(B,\ul{\Hom}(Gk,A),A,A,\tau,\mu)$ is strict.
\end{itemize}
\item \begin{itemize}
\item $B=T$;
\item $(F,G)$ is an equivalence of categories.
\end{itemize}
\end{enumerate}
\end{theorem}

\end{document}